\documentclass{amsart}
\usepackage{amssymb} 
\usepackage{amsmath} 
\usepackage{amscd}
\usepackage{amsbsy}
\usepackage{comment}
\usepackage[matrix,arrow]{xy}

\DeclareMathOperator{\End}{End}
\DeclareMathOperator{\rk}{rank}

\DeclareMathOperator{\NS}{NS}

\DeclareMathOperator{\GL}{GL}
\DeclareMathOperator{\SL}{SL}
\DeclareMathOperator{\Sel}{Sel}

\DeclareMathOperator{\rank}{rank}
\DeclareMathOperator{\Pic}{Pic}

\DeclareMathOperator{\Gal}{Gal}

\newcommand{\Q}{{\mathbb Q}}
\newcommand{\Z}{{\mathbb Z}}
\newcommand{\cA}{\mathcal{A}}

\newcommand{\OO}{{\mathcal O}}

\begin {document}

\newtheorem{thm}{Theorem}
\newtheorem{lem}{Lemma}[section]

\newtheorem{cor}[lem]{Corollary}

\theoremstyle{definition}

\theoremstyle{remark}

\title[]{Quadratic Chabauty for Modular Curves}
\author{Samir Siksek}
\address{Mathematics Institute\\
	University of Warwick\\
	Coventry\\
	CV4 7AL \\
	United Kingdom}

\email{s.siksek@warwick.ac.uk}

\date{\today}
\thanks{
The author is supported by
the EPSRC {\em LMF: L-Functions and Modular Forms} Programme Grant
EP/K034383/1.
}

\subjclass[2000]{Primary 11G30}

\begin {abstract}
Let $X/\Q$ be a curve of genus $g \ge 2$ with Jacobian $J$ and let $\ell$ be a prime
of good reduction. Using Selmer varieties, Kim defines a decreasing sequence
\[
X(\Q_\ell) \supseteq X(\Q_\ell)_1 \supseteq X(\Q_\ell)_2 \supseteq \cdots 
\]
all containing $X(\Q)$.
Thanks to the work of Coleman, the \lq Chabauty set\rq\ $X(\Q_\ell)_1$ is known
to be finite provided the \lq Chabauty condition\rq\ $\rank J(\Q) < g$ holds.
In this case one has a practical strategy that often succeeds in computing the set of rational points $X(\Q)$.
Balakrishnan and Dogra have recently shown that the \lq quadratic Chabauty set\rq\ $X(\Q_\ell)_2$
is finite provided 
\[
\rank J(\Q) < g + \rank \NS(J)-1,
\]
where $\NS(J)$ is the 
N\'eron-Severi group of $J/\Q$. In view of this it is interesting to give families of 
curves where $\rank \NS(J) \ge 2$ and where therefore quadratic Chabauty is
more likely to succeed than classical Chabauty. In this note we show that this is indeed the 
case for all modular curves of genus $\ge 3$.
\end {abstract}
\maketitle

\section{Introduction}

Let $X$ be a smooth algebraic curve of genus $g \ge 2$ defined over $\Q$ and write
$J$ for its Jacobian. 
Let $\ell$ be a 
prime of good reduction for $X$. The method of Chabauty and Coleman \cite{Coleman} yields
a subset $X(\Q_\ell)_1 \subseteq X(\Q_\ell)$ (which we may term the \lq Chabauty set\rq)  containing $X(\Q)$
that can be explicitly described in terms of Coleman integrals.
Coleman \cite{Coleman} showed that $X(\Q_\ell)_1$ is finite provided
the \lq Chabauty condition\rq\ holds:
\begin{equation}\label{eqn:chab1}
\rank J(\Q) < g .
\end{equation}
This implies that $X(\Q)$ is finite. Of course
Faltings' Theorem \cite{Fa1} asserts the
finiteness of $X(\Q)$ without assumptions beyond $g\ge 2$. 
However, Faltings' Theorem is ineffective. In contrast (see for example
\cite{MP} or \cite{Sikiii})  
the method of Chabauty and Coleman often allows for the computation of $X(\Q)$ 
provided the Chabauty condition \eqref{eqn:chab1} is satisfied
and one knows a Mordell--Weil basis for $J(\Q)$ (or even for a subgroup of
$J(\Q)$ of finite index).

In \cite{Kim}, Kim defined a family of Selmer varieties $\Sel(U_n)$ giving a decreasing sequence of subsets
\[
X(\Q_\ell) \supseteq X(\Q_\ell)_1 \supseteq X(\Q_\ell)_2 \supseteq \cdots \supseteq X(\Q),
\]
which can be computed in terms of iterated Coleman integrals. This offers hope 
of a strategy to compute $X(\Q)$ even if the Chabauty condition \eqref{eqn:chab1} fails.
Indeed, the conjectures of Bloch
and Kato imply that $X(\Q_\ell)_n$ is finite for $n$ sufficiently large \cite{Kim}.
The following theorem of Balakrishnan and Dogra \cite[Lemma 3]{BalD}
gives a  criterion for the finiteness of $X(\Q_\ell)_2$. 
\begin{thm}[Balakrishnan \& Dogra \cite{BalD}]
Suppose
\begin{equation}\label{eqn:chab2}
\rk J(\Q) < g+ \rk \NS(J)-1,
\end{equation}
where $\NS(J)$ is the N\'eron--Severi rank of $J$ over $\Q$. 
 Then $X(\Q_\ell)_2$ is finite.
\end{thm}
Indeed, Balakrishnan \& Dogra 
explain how, under condition \eqref{eqn:chab2}, a function 
can be constructed that 
cuts out $X(\Q_\ell)_2$
explicitly (following \cite{BBM} this method is termed \emph{Quadratic Chabauty}). They then apply this method to determine the rational
points on several genus $2$ curves
that have the form
\begin{equation}\label{eqn:Dogra}
X \; : \; y^2=a_6 x^6+ a_4 x^4+ a_2 x^2+a_0,
\end{equation}
where $\rk J(\Q)=2$. 

The method of Balakrishnan \& Dogra at first seems rather special.
For a \lq generic\rq\ curve $X$ it is known that the N\'eron--Severi
group has rank $1$ (and is in fact spanned by the class of the 
principal polarization). Hence the inequality
\eqref{eqn:chab2} is equivalent to the inequality \eqref{eqn:chab1},
and the method of Balakrishnan \& Dogra does not yield any more than 
the classical
Chabauty. Arguably modular curves are the most interesting family of 
curves. In this note, we give an explicit sufficiency criterion, based on 
the theorem of Balakrishnan \& Dogra, for the finiteness of $X(\Q_\ell)_2$
when $X$ is a modular curve. In fact we show that for any modular curve
$X/\Q$ of genus $\ge 3$ we have $\rank \NS(J) \ge 2$, thus
quadratic Chabauty is more likely to succeed than classical Chabauty.

\section{The Balakrishnan--Dogra Criterion and Modular Curves}

Let $G$ be a congruence subgroup of $\SL_2(\Z)$ and let $X=X_G$ be
the corresponding modular curve which we shall suppose is defined
over $\Q$ (sufficient conditions on $G$ for this to hold are given 
in \cite[Section 1.2]{Rohrlich}). Assume that $X$ has genus $g \ge 2$ and
write $J$ for the Jacobian of $X$.
Let $f_1,\dotsc,f_n$ be representatives
of the Galois orbits of weight $2$ cuspidal eigenforms for $G$ (these may
be computed via the modular symbols algorithm as in \cite{Stein}). 
The theory of Eichler--Shimura attaches to each
$f_i$ an abelian variety of $\GL_2$-type, which we shall denote
by $\cA_i$. 
Then 
\begin{equation}\label{eqn:factor}
J \approx \cA_{1} \times \cdots \times \cA_{n},
\end{equation}
where $\approx$ denotes isogeny over $\Q$.
As $\cA_i$ is of $\GL_2$-type, we have $\End(\cA_i) \otimes \Q$
 is a number field
$F_i$ of degree $\dim \cA_i$, which is either totally real
or CM. The fields $F_i$ are generated by the Hecke eigenvalues of $f_i$.
We arrange the $f_i$ so that $F_1,\dotsc,F_m$ are totally
real and $F_{m+1},\dots,F_n$ are CM.  
\begin{cor}
With the above notation and assumptions,
suppose 
\begin{equation}\label{eqn:betterchab1}
\rank J(\Q)< -1+ 2 \sum_{i=1}^m [F_i:\Q]+\frac{3}{2}\sum_{i=m+1}^n [F_i:\Q].
\end{equation}
Then, for any prime $\ell$ of good reduction for $X$, 
the set $X(\Q_\ell)_2$ is finite.
\end{cor}
\begin{proof}
Note that
\[
g=\dim(J)=\sum_{i=1}^n \dim \cA_i=\sum_{i=1}^n [F_i:\Q].
\]
Thus \eqref{eqn:betterchab1} maybe rewritten as
\begin{equation}\label{eqn:betterchab}
\rank J(\Q)< (g-1)+\sum_{i=1}^m [F_i:\Q]+\frac{1}{2}\sum_{i=m+1}^n [F_i:\Q].
\end{equation}
It is sufficient to show, thanks to the theorem of Balakrishnan \& Dogra,
that
\[
\rank \NS(J) \ge \sum_{i=1}^m [F_i:\Q]+\frac{1}{2}\sum_{i=m+1}^n [F_i:\Q].
\]
We now invoke the isomorphism of $\Q$-vector spaces
\[
\NS(J) \otimes \Q \cong \End^{(s)}(J) \otimes \Q,
\]
where $\End^{(s)}(J)$ is the subring of $\End(J)$ left invariant by
the Rosati involution (see \cite[Section 1]{Giovanetti} for example). 
However, the Rosati involution restricted to $\End(\cA_i)$ induces
complex conjugation on $F_i$ (\cite[page 196]{Pyle}). 
Thus we see that $\End^{(s)}(J) \otimes \Q$
contains the algebra 
\begin{equation}\label{eqn:algebra}
F_1\times \cdots \times F_m \times E_{m+1} \times \cdots \times E_n,
\end{equation}
where, for $m+1 \le i \le n$, we denote the maximal totally real
subfield of $F_i$ by $E_i$. As $F_i$ is CM, we have $[E_i:\Q]=[F_i:\Q]/2$.
The corollary follows.
\end{proof}
It is evident that \eqref{eqn:betterchab} (and hence \eqref{eqn:betterchab1})
is weaker than
Chabauty's condition \eqref{eqn:chab1} whenever $g \ge 3$. If $g=2$
then \eqref{eqn:betterchab} 
is weaker if either $J$ has real multiplication, or is isogenous to 
a product of two elliptic curves.

\medskip

\noindent \textbf{Remark}. If $J$ has multiplicity $1$ (meaning
that the
simple factors $\cA_i$ are pairwise non-isogenous) then 
$\End(J) \cong \End(\cA_1) \times \cdots \times \End(\cA_n)$
and the 
above argument in fact shows that
$\NS(J) \otimes \Q$ is isomorphic to the algebra \eqref{eqn:algebra}.

\end{document}